\documentclass[11pt]{amsart}

\usepackage[english]{babel}
\usepackage{amsmath, amsfonts, amssymb, amsthm, epsfig, graphicx, url, hyperref, color, tikz, cancel}

\newtheorem{theorem}{Theorem}[section]

\newtheorem{definition}{Definition}[section]

\newtheorem{remark}{Remark}[section]

\newcounter{theor}

\newtheorem{thm}[theor]{Theorem}
\newtheorem{cor}[theor]{Corollary}

\newcounter{lemm}

\newtheorem{lem}[lemm]{Lemma}

\DeclareMathOperator{\inter}{int}
\DeclareMathOperator{\relinter}{relint}

\def\conv{\mathop\mathrm{conv}\nolimits}

\def\supp{\mathop\mathrm{supp}\nolimits}

\def\K{\mathcal{K}}
\def\w{\mathcal{W}}

\def\R{\mathbb{R}}

\def\Z{\mathbb{Z}}

\def\V{\mathrm{V}}
\def\vol{\mathrm{vol}}
\def\F{\mathcal{F}}

\def\W{\mathrm{W}}

\def\e{\mathrm{e}}
\newcommand{\dlat}{\mathrm{d}}

\def\e{\mathrm{e}}

\numberwithin{equation}{section}

\begin{document}
\title[$L_p$ Brunn-Minkowski type inequalities for a general class of functionals]{On $L_p$ Brunn-Minkowski type inequalities for a general class of functionals}

\author[L. Gordo Malag\'on]{Lidia Gordo Malag\'on}
\address{Departamento de Matem\'aticas, Universidad de Murcia, Campus de
Espinar\-do, 30100 Murcia, Spain}
\email{lidia.gordom@um.es} \email{jesus.yepes@um.es}

\author[J. Yepes Nicol\'as]{Jes\'us Yepes Nicol\'as}

\thanks{Both authors are supported by the grant PID2021-124157NB-I00, 
funded by MCIN/AEI/10.13039/501100011033/ ``ERDF A way of making Europe''. The second author is also supported by the grant  ``Proyecto financiado por la CARM a través de la convocatoria de Ayudas a proyectos para el desarrollo de investigación científica 
y técnica por grupos competitivos, incluida en el Programa Regional de Fomento de la Investigación Científica y Técnica (Plan de Actuación 2022) de la Fundación Séneca-Agencia de Ciencia y Tecnología de la Región de Murcia, REF. 21899/PI/22'', as well as by the grant RYC2021-034858-I funded by MCIN/AEI/10.13039/501100011033 and by the ``European Union NextGenerationEU/PRTR''}



\begin{abstract}
In this work, the $L_p$ version (for $p> 1$) of the dimensional Brunn-Minkowski inequality for the standard Gaussian measure $\gamma_n(\cdot)$ on $\mathbb{R}^n$ is shown. More precisely, we prove that for any $0$-symmetric convex sets with nonempty interior, any $p>1$, and every $\lambda \in (0,1)$,  
\[
\gamma_n\bigl((1-\lambda)\cdot K+_p \lambda \cdot L\bigr)^{p/n} \geqslant (1-\lambda ) \gamma_n(K)^{p/n} + \lambda \gamma_n(L)^{p/n},
\]
with equality, for some $\lambda \in (0,1)$ and $p>1$, if and only if $K=L$. This result, recently established without the equality conditions by Hosle, Kolesnikov and Livshyts, by using a different and functional approach, turns out to be the $L_p$ extension of a celebrated result for the Minkowski sum (that is, for $p=1$) by Eskenazis and Moschidis (2021) on a problem by Gardner and Zvavitch (2010).

Moreover, an $L_p$ Brunn-Minkowski type inequality is obtained for the classical Wills functional $\mathcal{W}(\cdot)$ of convex bodies.

These results are derived as a consequence of a more general approach, which provides us with other remarkable examples of functionals satisfying $L_p$ Brunn-Minkowski type inequalities, such as different absolutely continuous measures with radially decreasing densities.
\end{abstract}

\maketitle

\section{Introduction}
We work in the $n$-dimensional Euclidean space $\R^n$ endowed with the standard scalar pro\-duct $\langle \cdot, \cdot \rangle$, and Euclidean norm $\|\cdot\|$. The unit (closed) ball is denoted by $B_n$, a set $A$ is referred to as $0$-symmetric
if $A=-A$, and a set of the form $\lambda A$, for some $\lambda \geqslant 0$, is called a dilatate of $A$. Let $\mathcal{K}^n$ be the set of all convex bodies, i.e., nonempty compact convex sets, in $\R^n$, and let $\mathcal{K}^n_0$ be the subfamily of $\mathcal{K}^n$ of all convex bodies containing the origin as an interior point.
The volume of a measurable set $A \subset \R^n$, i.e., its $n$-dimensional Lebesgue measure, is denoted by $\vol(A)$ (when integrating, as usual, $\dlat x$ stands for $\dlat \vol(x)$) and, in particular, we write $\kappa_n=\vol(B_n)$. Moreover, with $\inter A$ and $\conv A$, we represent the interior and the convex hull of $A$, respectively, and with $\relinter A$ we denote its relative interior, i.e., the interior of $A$ relative to its affine hull. 
The orthogonal projection of $A$ onto a (vector) subspace $H$ is denoted by $A|H$ and with $H^{\perp}$ we represent the orthogonal complement of $H$.
Besides, if $K \in \K^n_0$ is a convex body containing the origin in its interior, then the polar body $K^*$ of $K$ is defined by $K^*=\{ x \in \R^n : \langle x, y \rangle \leqslant 1 $ for all $y \in K\}$. Finally, $K+L=\{x+y: x \in K, \,y \in L\}$ is the Minkowski sum of any nonempty sets $K, L \subset\R^n$.

Relating the volume $\vol(\cdot)$ with the Minkowski addition of two convex bodies $K,L\in\K^n$, one is led to the well-known \emph{Brunn-Minkowski inequality}, which is one of the cornerstones of the Brunn-Minkowski theory. It assures that, for any $\lambda\in (0,1)$, 
\begin{equation}\label{e: BM}
{\vol\bigl((1-\lambda)K+ \lambda L\bigr)^{{1}/{n}} \geqslant (1-\lambda ) \vol(K)^{{1}/{n}} + \lambda \vol(L)^{{1}/{n}}}.
\end{equation}
Equality for some $\lambda \in (0,1)$ holds if and only if $K$ and $L$ either lie in parallel hyperplanes or are homothetic.

In $1962$, Firey \cite{F2} introduced the following concept of $p$\emph{-sum} or $L_p$ \emph{addition}: for two convex bodies containing the origin $K$, $L \subset \R^n$ and $1\leqslant p \leqslant \infty$ fixed there exists a unique convex body $K +_p L $ whose support function is given by 
\begin{equation}\label{e: p-sum Firey}
{h(K+_p L,\cdot)=\bigl(h(K,\cdot)^p+h(L,\cdot)^p\bigr)^{1/p}}.
\end{equation}
When $p=\infty$ this must be considered as its limit case, i.e., $h(K +_{\infty} L, \cdot) = \max\bigl\{h(K, \cdot), h(L, \cdot)\bigr\}$, as usual.
We recall that $h(K,u)=\max\bigl\{ \langle x,u \rangle :x \in K \bigr\}$ for all $u \in \R^n$ is the support function of a convex body $K\subset\R^n$. And for convenience, in this context, we define the $p$-scalar multiplication by $\lambda \cdot K=\lambda^{1/p}K$, for any $\lambda\geqslant 0$.  

Clearly, $+_1$ is the standard Minkowski addition, whereas $+_\infty$ yields
\begin{equation*}
K +_{\infty} L = \conv (K \cup L). 
\end{equation*}

Noting that the support function of a set equals the support function of its convex hull, we point out that Firey's definition of $p$-sum \eqref{e: p-sum Firey} requires assuming both convexity, since it is given in terms of the support functions of the convex bodies involved, and that the sets contain the origin, in contrast to what happens for the classical Minkowski addition. So, in $2012$, Lutwak, Yang and Zhang \cite{LYZ} extended the above notion of $p$-sum to the case of arbitrary subsets of $\R^n$. Then, if $K, L \subset \R^n$ are nonempty sets and $1\leqslant p < \infty$, they defined the $p$-sum by
\begin{equation*}\label{e: p-sum LYZ}
{K +_p L=\left\{(1-\mu)^{1/q}x+\mu^{1/q}y: x\in K,y\in L, \, \mu\in[0,1]\right\}},
\end{equation*}
where $q$ denotes the H\"older conjugate of $p$, i.e., such that 
\begin{equation}\label{e: conjugate Holder}
\frac{1}{p}+\frac{1}{q}=1.
\end{equation}
In \cite{LYZ} it is shown that the latter definition coincides with the one given by Firey when $K$ and $L$ are convex bodies containing the origin. Furthermore, when $p=1$ (and hence $q=\infty$), the coefficients $(1-\mu)^{1 / q}$  and $ \mu^{1 / q}$ must be interpreted as $1$ for all $0 \leqslant \mu \leqslant 1$. Finally, following the approach taken in \cite{LYZ}, we omit the case when $p=\infty$ (and so $q=1$), since for such a value of $p$ all the results hold trivially. Thus, throughout the rest of the manuscript, whenever $p \geqslant 1$ is mentioned, we will refer to a real number $p \geqslant 1$.

The $L_p$ version of the Brunn-Minkowski inequality \eqref{e: BM} was originally established by Firey \cite{F2} for convex bodies containing the origin, and later extended by Lutwak, Yang and Zhang (see \cite[Theorem 4]{LYZ}) to arbitrary nonempty compact sets. It states the following:
\begin{thm}
Let $K, L \subset\R^n$ be nonempty compact sets and $p > 1$. Then, for all $\lambda\in(0,1)$,
\begin{equation*}
    {\vol\bigl((1-\lambda){\cdot} K +_p \lambda{\cdot}L\bigr)^{p/n} \geqslant(1-\lambda)\vol(K)^{{p}/n}+\lambda\vol(L)^{{p}/n}}.
\end{equation*}
\end{thm}

Approximately three decades after Firey introduced the concept of $p$-sum for convex bo\-dies (containing the origin), Lutwak \cite{L1, L2} developed a fruitful theory, the so-called $L_p$ Brunn-Minkowski theory, which nowadays is a very active area of research and the starting point for new developments and generalizations. For more information on the $L_p$ Brunn-Minkowski theory and its consequences, we refer the reader to \cite[Section 9.1]{Sch} and the references therein.

\smallskip

This paper mainly focuses on finding $L_p$ versions of different Brunn-Minkowski type inequalities. In particular, we aim to extend to the $L_p$ setting recent Brunn-Minkowski inequalities for both the \emph{Gaussian measure} (and, more generally, some absolutely continuous measures) and the classical \emph{Wills functional}, among other functionals.

To this end, first we deal with the case of the standard Gaussian measure $\gamma_n$ on $\R^n$, given by 
\begin{equation*}
{\mathrm{d}\gamma_n(x)=\frac{1}{(2\pi)^{n/2}}e ^\frac{-\|x\|^2}{2} \;\mathrm{d}x},
\end{equation*}
for which a dimensional Brunn-Minkowski type inequality has been explored. More precisely, taking into account \eqref{e: BM}, the following natural question arises: does the inequality
\begin{equation}\label{e: Gaussian BM}
\gamma_n\bigl((1-\lambda)K+\lambda L\bigr)^{1/n} \geqslant (1-\lambda) \gamma_n(K)^{1/n} + \lambda \gamma_n(L)^{1/n} 
\end{equation}
hold for any closed convex sets $K,L\subset\R^n$ containing
the origin and all $\lambda\in(0,1)$? In \cite{GZ} Gardner and Zvavitch first observed that \eqref{e: Gaussian BM} is not true without any restriction on the position of the sets. To see this it is enough to consider $K=B_n$ and $L=B_n+x$, with $x$ large enough. But, if the sets $K$ and $L$ contain the origin, Gardner and Zvavitch conjectured that this Gaussian Brunn-Minkowski inequality would be true. In particular, they proved that this inequality holds true for dimension 1, for coordinate boxes containing the origin, that is, boxes (containing the origin) whose sides are parallel to the coordinate axes, and when either $K$ or $L$ is a slab containing the origin. Although, in 2013, Nayar and Tkocz \cite{NT} showed that this conjecture is in general false, the possibility of such an inequality being true for all $0$-symmetric convex sets remained open until 2021, when Eskenazis and Moschidis \cite{EM} proved the following celebrated result:

\begin{thm}\label{t: EskMosch Gaussian BM}
Let $K, L \subset \mathbb{R}^n$ be $0$-symmetric closed convex sets with nonempty interior. Then, for all $\lambda \in (0,1)$,  
\begin{equation*} 
\gamma_n\bigl((1-\lambda) K+ \lambda  L\bigr)^{1/n} \geqslant (1-\lambda ) \gamma_n(K)^{1/n} + \lambda \gamma_n(L)^{1/n}.
\end{equation*}
Equality, for some $\lambda \in (0,1)$, holds if and only if $K=L$.
\end{thm}

In relation to this result, here we study its corresponding $L_p$ version, for $p > 1$, obtained as a consequence of our main result, Theorem \ref{teoremafunctionalFnonconvex}, in the setting of arbitrary general functionals defined on a family of subsets of $\R^n$. In this regard, we show the following $L_p$ Gaussian Brunn-Minkowski inequality, which was established by Hosle, Kolesnikov and Livshyts \cite{HKL} in 2020 (without the equality conditions), by using a different approach. There they considered and studied local versions of related functional inequalities. Their elegant method, which is often referred to in the literature as a local-to-global principle, was later exploited by Eskenazis and Moschidis in \cite{EM}.

\begin{theorem}\label{t:Gaussian measure symmetric convex Lp BM}
Let $K, L \subset \mathbb{R}^n$ be $0$-symmetric closed convex sets with nonempty interior and $p > 1$. Then, for all $\lambda \in (0,1)$,  
\begin{equation}\label{e:Gaussian measure symmetric convex Lp BM} 
\gamma_n\bigl((1-\lambda)\cdot K+_p \lambda \cdot L\bigr)^{p/n} \geqslant (1-\lambda ) \gamma_n(K)^{p/n} + \lambda \gamma_n(L)^{p/n}.
\end{equation}
Equality, for some $\lambda \in (0,1)$ and $p>1$, holds if and only if $K=L$.
\end{theorem}

Moreover, in Section \ref{s: L_p BM Gaussian}, we show that the above $L_p$ Gaussian Brunn-Minkowski inequality \eqref{e:Gaussian measure symmetric convex Lp BM} is true when $K$ and $L$ are either \emph{weakly unconditional} measurable sets or convex bodies with many symmetries (see Theorems \ref{t: Gaussian Lp BM weakly} and \ref{t: Gaussian Lp BM BK}, and Definition \ref{def: weakly sets}).

Furthermore, when computing the volume of the Minkowski sum $K + \lambda B_n$, with $\lambda \geqslant 0$, one is led to the classical \emph{Steiner formula} of $K$, namely, 
\begin{equation}\label{eq:Steiner formula}
    \vol(K+\lambda B_n)= \sum_{i=0}^n\binom{n}{i} \W_{i}(K) \lambda^i.
\end{equation}
The coefficients $\W_i(K)$ are the \emph{quermassintegrals} of $K$, and they are a special case of the more general defined \emph{mixed volumes}, for which we refer to \cite[Section 5.1]{Sch}. In particular $\W_0(K) = \vol(K)$, $\W_n(K) = \vol(B_n)$, $n \W_1(K) = \mathrm{S}(K)$ is the surface area of $K$ and $(2/\kappa_n)\W_{n-1}(K) = \mathrm{b}(K)$ is the mean width of $K$ (see \cite[p. 50]{Sch}). 
 
In \cite{Mc}, McMullen introduced the normalized quermassintegrals, defined as
\begin{equation*}
\V_{n-i}(K)=\binom{n}{i}\frac{\W_{i}(K)}{\kappa_i},
\end{equation*}
and suggested referring to these measures as the \emph{intrinsic volumes} of $K$, since if $K$ has dimension $k$ then $\V_k(K)$ coincides with the $k$-dimensional volume of $K$, and moreover the intrinsic volumes are independent of the dimension of the space in which $K$ is embedded (see e.g. \cite[Section 6.4]{Gr}).

In $1973$ (see \cite{W}) it was introduced and studied the functional (nowadays usually referred to in the literature as the Wills functional) given by the sum of all the intrinsic volumes, i.e.,
\begin{equation*}\label{e: def Wills functional}
{\w(K)=\sum_{i=0}^{n}\V_i(K)},
\end{equation*}
because of its possible connection with the so-called lattice-point enumerator $\mathrm{G}(K) = \#(K \cap \Z^n)$, where $\#$ denotes the cardinality. In recent years, numerous interesting properties of this functional have been investigated (see, for instance, \cite{H, H2, W, W2, W3}).

Taking into account that all the summands of the Wills functional are $(1/n)$\emph{-concave} (which means that they satisfy \eqref{e: BM}), since one has (see e.g. \cite[Theorem 7.4.5]{Sch}) that 
\begin{equation*}
    \V_i \bigl( (1-\lambda) K +\lambda L \bigr)^{1/i} \geqslant (1-\lambda) \V_i(K)^{1/i} + \lambda \V_i(L)^{1/i}
\end{equation*}
for any $K, L \in \K^n$ and all $\lambda \in (0,1)$, for any $i=1, \dots, n$, it is natural to wonder whether a dimensional Brunn-Minkowski inequality also holds for the Wills functional. Although, as seen in \cite{AHY}, the Wills functional is unfortunately not $(1/n)$-concave, in $2021$ Alonso-Gutiérrez, Hernández Cifre and the second-named author \cite{AHY} proved that an alternative inequality holds when a constant (depending on the dimension) is introduced:
\begin{thm} \label{t: BM Wills}
Let $K, L \in \K^n$ be convex bodies. Then, for all $\lambda \in (0,1)$,
\begin{equation*}
\w\bigl((1-\lambda)K+\lambda L\bigr)^{1/n}
\geqslant{\frac{{1}}{{(n!)^{1/n}}}}{\left((1-\lambda)\w(K)^{1/n}
+\lambda\w(L)^{1/n}\right)}.
\end{equation*}
\end{thm}

Here we show that the above result admits an extension for the $p$-sum, which yields the corres\-ponding $L_p$ Brunn-Minkowski type inequality. More precisely, we prove the following:

\begin{theorem}\label{t:classic Wills functional Lp BM}
Let $K,L \in \K^n$ be convex bodies and $p >1$. Then, for all $\lambda \in (0,1)$, 
\begin{equation*}
{\w \bigl((1-\lambda)\cdot K+_p \lambda \cdot L\bigr)^{{p}/{n}} \geqslant \frac{1}{(n!)^{{p}/{n}}} \left((1-\lambda ) \w(K)^{{p}/{n}} + \lambda \w(L)^{{p}/{n}}\right)}.
\end{equation*}
\end{theorem}

We would like to point out that Theorems \ref{t:Gaussian measure symmetric convex Lp BM} and \ref{t:classic Wills functional Lp BM} are derived as a consequence of a more general approach, shown in Section \ref{s: General L_p BM}. There we deal with a functional (defined on a family of subsets of $\R^n$) that is both increasing and \emph{sub-homogeneous} (see Section \ref{s: General L_p BM} for the definition), and then we show that if a certain pair of subsets satisfies a Brunn-Minkowski type inequality then the corresponding $L_p$ version of it also holds (see Theorem \ref{teoremafunctionalFnonconvex}, where the precise equality conditions when $p>1$ are further obtained). 
At this point, we would like to notice that such a sub-homogeneity property is often fulfilled by different functionals. For example, as it will be discussed in Section \ref{s: L_p BM Gaussian}, we may consider any absolutely continuous measure on $\R^n$ with \emph{radially decreasing} density function, since such a measure is sub-homogeneous of degree $n$. Thus, we have a rich sample of absolutely continuous measures being sub-homogeneous of degree $n$, in contrast with the (degree $n$) homogeneous case, where the Lebesgue measure is the sole example of such a measure with continuous density function. In this regard, it is fair to mention that the homogeneous case of ($L_p$ Brunn-Minkowski inequalities for) general functionals was previously obtained in \cite{ZX}, where the authors provide some examples of applications of their result.

\smallskip

The paper is organized as follows: in Section \ref{s: General L_p BM} we prove our main result, for general monotonous and sub-homogeneous functionals defined on a family of subsets of $\R^n$, and a couple of sets satisfying a Brunn-Minkowski inequality for the Minkowski sum, by showing that then the corresponding Brunn-Minkowski inequality for the $L_p$ addition holds for those sets. In Section \ref{s: L_p BM Gaussian} 
we discuss various different consequences of such a result for general absolutely continuous measures with radially decreasing densities, and in particular for the Gaussian measure (obtaining among others Theorem \ref{t:Gaussian measure symmetric convex Lp BM}), whereas the corresponding $L_p$ Brunn-Minkowski type inequality for the (generalized) Wills functional is derived in Section \ref{s: BM Wills functional}.

\smallskip

\section{A general $L_p$ Brunn-Minkowski inequality}\label{s: General L_p BM}

In this section we prove our main result, which provides an $L_p$ Brunn-Minkowski inequality for a certain class of functionals (defined on a family of nonempty subsets of $\R^n$) and a pair of sets for which one has a Brunn-Minkowski inequality for the Minkowski addition. In order to state our result in a precise way, we need first the following definition.

\begin{definition}
We say that a nonnegative functional $\F: \mathcal{A}  \longrightarrow \R_{\geqslant0}$, defined on a family $\mathcal{A}$ of nonempty subsets of $\R^n$ closed under dilatations, is sub-homogeneous of degree $1/\alpha$, $\alpha\neq 0$, if 
\begin{equation}\label{e: sub-homog}
\F(r K) \leqslant r ^{1/\alpha} \F(K)
\end{equation} 
whenever $r \geqslant 1$. Moreover, it is strictly sub-homogeneous of degree $1/\alpha$ if \eqref{e: sub-homog} is strict provided that $\F(K)>0$ and $r > 1$. 

Analogously, we say that $\F: \mathcal{A}  \longrightarrow \R_{\geqslant0}$ is super-homogeneous of degree $1/\alpha$, $\alpha\neq0$, if 
\begin{equation}\label{e: super-homoge}
\F(r K) \geqslant r ^{1/\alpha} \F(K)
\end{equation} 
whenever $r \geqslant 1$. Furthermore, it is strictly super-homogeneous of degree $1/\alpha$ if \eqref{e: super-homoge} is strict provided that $\F(K)>0$ and $r > 1$. 
\end{definition}

Before stating and proving our result, we need to recall the following lemma, shown in \cite{LYZ}, which relates the $L_p$ and classical convex combinations in terms of set inclusion.
\begin{lem}\label{lemaLYZ}
Let $K,L \subset \R^n$ be nonempty sets and $p>1$. Then, for all $\lambda \in (0,1)$,
\begin{equation*}
(1 - \lambda) K + \lambda L \subset (1 - \lambda) \cdot K +_p \lambda \cdot L.
\end{equation*}
Equality, for some $\lambda \in (0,1)$, implies that $\conv K=\conv L$ and $0\in \conv K \cap \conv L$, provided that $K$ and $L$ are compact sets. 
\end{lem}

We are now in a position to show our main result, for which we exploit the original approach followed by Lutwak, Yang and Zhang in the proof of \cite[Theorem 4]{LYZ}.
\begin{theorem} \label{teoremafunctionalFnonconvex}
Let $\mathcal{A}$ be a family of nonempty subsets of $\R^n$ closed under linear combinations, and let $\F: \mathcal{A}  \longrightarrow \R_{\geqslant0} $ be a nonnegative functional that is increasing under set inclusion and sub-homogeneous of degree $1/\alpha$ for some $\alpha>0$. 
Let $K,L\in\mathcal{A}$ with $\F(K)\F(L)>0$ be such that, for all $\lambda \in (0,1)$,
\begin{equation}\label{e:BM functional F} 
\F \bigl((1-\lambda) K+\lambda L \bigr) \geqslant C \bigl((1-\lambda) \F(K)^\alpha + \lambda \F(L)^\alpha \bigr)^{1/\alpha}
\end{equation}
for some constant $C>0$.
Then, for any $p\geqslant 1$ and all $\lambda \in (0,1)$, we have
\begin{equation}\label{e: LpBMforFnonconvex}
\mathcal{F}\bigl((1-\lambda)\cdot K +_p \lambda \cdot L\bigr) \geqslant C \bigl( (1-\lambda) \mathcal{F}(K)^{{p \alpha}} + \lambda \mathcal{F}(L)^{{p \alpha}}\bigr)^{1/p \alpha} 
\end{equation} 
whenever $(1-\lambda)\cdot K+_p\lambda\cdot L \in \mathcal{A}$. 

\smallskip

\noindent
Equality, for some $\lambda \in (0,1)$ and $p>1$, implies that $K$ and $L$ satisfy \eqref{e:BM functional F} with equality for certain $\bar{\lambda}\in(0,1)$.
Moreover, if $\F$ is strictly increasing and $K$ and $L$ are compact then $\conv K$ and $\conv L$ contain the origin and are dilatates of each other. If further $\F$ is strictly sub-homogeneous of degree $1/\alpha$ then $\conv K = \conv L$.
\end{theorem}

\begin{remark}
For the case where $K, L \in \K^n$ are convex bodies containing the origin, and $C=1$, equality in \eqref{e: LpBMforFnonconvex} holds if and only if $K=L$. 
\end{remark}

\begin{proof} 
Set $p>1$ and let $\lambda,\mu\in(0,1)$.
From the monotonicity of the functional $\F$ 
jointly with Lemma \ref{lemaLYZ} applied to the sets $(1-\mu)^{-1/p}(1-\lambda)\cdot K$ and $\mu^{-1/p} \lambda \cdot L$, we get
\begin{equation}\label{e:F_property i + LEMMA LYZ}
\begin{split}
\F \bigl(\left(1-\lambda \right) \cdot K +_p \lambda \cdot L\bigr) &=
\F\bigl((1-\mu)^{1/p}(1-\mu)^{-1/p} (1-\lambda)\cdot K  +_p \mu^{1/p} \mu^{-1/p} \lambda \cdot L \bigr)\\  
&\geqslant \F\bigl((1-\mu)(1-\mu)^{-1/p} (1-\lambda)\cdot K  + \mu \mu^{-1/p} \lambda \cdot L \bigr)\\ 
&=\F\bigl((1-\mu)^{1/q} (1-\lambda)^{1/p} K  + \mu^{1/q} \lambda^{1/p} L \bigr).
\end{split}
\end{equation}
Now, setting $t=(1-\mu)^{1/q}(1-\lambda)^{1/p}$ and $s=\mu^{1/q}\lambda^{1/p}$, we have that $t+s \leqslant 1$ by H\"older's inequality and then from the sub-homogeneity of 
$\F$, applied with $r=1/(t+s)$,
we obtain
\begin{equation}\label{e:subhomog F}
\F(tK+sL)\geqslant (t+s)^{1/\alpha}\,
\F\left( \frac{t}{t+s}K + \frac{s}{t+s}L\right).
\end{equation}
Then, from \eqref{e:BM functional F} with $\bar{\lambda}=s/(t+s)$ 
we get
\begin{equation}\label{e:BM F}
(t+s)^{1/\alpha}\F\left( \frac{t}{t+s} K + \frac{s}{t+s}L\right) 
\geqslant C\bigl(t \F(K)^{\alpha} + s \F(L)^\alpha \bigr)^{1/ \alpha}, 
\end{equation}
and thus, from \eqref{e:F_property i + LEMMA LYZ}, \eqref{e:subhomog F} and \eqref{e:BM F}, we have that
\begin{equation}\label{e:Fadditivepreciset,snonconvex}
 \mathcal{F} \bigl( (1-\lambda)\cdot K +_p \lambda \cdot L \bigr) \geqslant C \left((1-\mu)^{1/q}(1-\lambda)^{1/p}\F(K)^{\alpha} + \mu^{1/q} \lambda^{1/p}\F(L)^{\alpha}\right)^{1/ \alpha}.
\end{equation}
Since the latter inequality holds for any $\mu\in (0,1)$, we may in particular take
\begin{equation*}\label{e:precise value mu}
\mu= \frac{\lambda \F(L)^{p\alpha}}{(1-\lambda) \F(K)^{p\alpha}+\lambda \F(L)^{p \alpha}},
\end{equation*}
and therefore, using that $q$ is the Hölder conjugate of $p$ (cf. \eqref{e: conjugate Holder}), \eqref{e:Fadditivepreciset,snonconvex} yields the desired inequality \eqref{e: LpBMforFnonconvex}.

Next, we study the equality conditions. So, if we assume equality in \eqref{e: LpBMforFnonconvex} for some $\lambda\in(0,1)$, then we have equality in \eqref{e:Fadditivepreciset,snonconvex} and thus in \eqref{e:F_property i + LEMMA LYZ}, \eqref{e:subhomog F} and \eqref{e:BM F}, and hence in particular $K$ and $L$ satisfy \eqref{e:BM functional F} with equality for some $\bar{\lambda}\in(0,1)$.
Moreover, from the strict monotonicity of $\F$, equality in \eqref{e:F_property i + LEMMA LYZ} implies that the sets 
$(1-\lambda)\cdot K  +_p \lambda \cdot L$ and 
$(1-\mu)^{1/q} (1-\lambda)\cdot K  + \mu^{1/q} \lambda \cdot L$
are the same. Thus, since $K$ and $L$ are compact, from the equality case of Lemma \ref{lemaLYZ} we obtain that 
\[
\conv \bigl((1-\mu)^{-1/p}(1-\lambda)\cdot K \bigr)= \conv \bigl(\mu^{-1/p} \lambda \cdot L\bigr)
\] 
and $0 \in \conv\bigl((1-\mu)^{-1/p}(1-\lambda)\cdot K \bigr) \cap \conv \bigl(\mu^{-1/p} \lambda \cdot L\bigr)$, 
that is, $\conv K$ and $\conv L$ are dilatates and contain the origin.

\smallskip

Furthermore, using the strict sub-homogeneity of $\F$ jointly with the fact that $\F(tK+sL)>0$, since $\F(K)\F(L)>0$ (cf. \eqref{e:subhomog F} and \eqref{e:BM F}), equality in \eqref{e:subhomog F} now implies that $t+s=1$. Hence, by the equality condition of Hölder's inequality, which was applied to the vectors 
\[
\bigl( (1-\mu)^{1/q}, \mu^{1/q}\bigr) \quad\text{and} \quad \bigl( (1-\lambda)^{1/p}, \lambda^{1/p}\bigr),
\] 
we have that $(1-\mu) = c (1-\lambda) $ and $\mu=c \lambda$ (for some $c\in\R$), which implies that $\mu=\lambda$.  
Thus, we get that 
\[
(1-\mu)^{-1/p}(1-\lambda)\cdot K=K
\quad \text{and} \quad 
\mu^{-1/p} \lambda \cdot L=L,
\]
and therefore $\conv K= \conv L$. This concludes the proof.
\end{proof}

\begin{remark}\label{r: A not necessarily closed under linear comb}
From the proof of the above result, we observe that one may consider an arbitrary family $\mathcal{A}$ of nonempty subsets of $\R^n$, i.e., not necessarily closed under linear combinations, provided that the pair of sets $K,L\in\mathcal{A}$ considered therein (in the conditions of the above result) are such that any linear combination $\lambda_1 K+\lambda_2 L\in\mathcal{A}$, for all $\lambda_1, \lambda_2\geqslant0$ with $\lambda_1+\lambda_2\leqslant 1$.

The same observation applies to the following result (Theorem \ref{teorema F alpha negative}).
\end{remark}

The Wills functional 
and the Gaussian measure are remarkable examples of functionals in the conditions of Theorem \ref{teoremafunctionalFnonconvex}, because both of them satisfy the properties therein for $\alpha = 1/n$. We will discuss about this in a more precise way in the forthcoming sections. 

\smallskip 

Now we state the corresponding inequality for the case when $\alpha < 0$. It is an analogue of the latter result, but for a functional $\mathcal{F}$ with the next properties: decreasing monotonicity, super-homogeneity and $\alpha$-convexity (cf. \eqref{e:BM functional F neg alpha}).
One may reproduce the proof of Theorem \ref{teoremafunctionalFnonconvex}, with the same steps, to show the following result.

\begin{theorem}\label{teorema F alpha negative}
Let $\mathcal{A}$ be a family of nonempty subsets of $\R^n$ closed under linear combinations, and let $\F: \mathcal{A}  \longrightarrow \R_{\geqslant0} $ be a nonnegative functional that is decreasing under set inclusion and super-homogeneous of degree $1/\alpha$ for some $\alpha<0$. 
Let $K,L\in\mathcal{A}$ with $\F(K)\F(L)>0$ be such that, for all $\lambda \in (0,1)$,
\begin{equation}\label{e:BM functional F neg alpha} 
\F \bigl((1-\lambda) K+\lambda L \bigr) \leqslant C \bigl((1-\lambda) \F(K)^\alpha + \lambda \F(L)^\alpha \bigr)^{1/\alpha}
\end{equation}
for some constant $C>0$.
Then, for any $p\geqslant 1$ and all $\lambda \in (0,1)$, we have
\begin{equation*}\label{e: LpBMforFnonconvex neg alpha}
\mathcal{F}\bigl((1-\lambda)\cdot K +_p \lambda \cdot L\bigr) \leqslant C \bigl( (1-\lambda) \mathcal{F}(K)^{{p \alpha}} + \lambda \mathcal{F}(L)^{{p \alpha}}\bigr)^{1/p \alpha} 
\end{equation*} 
whenever $(1-\lambda)\cdot K+_p\lambda\cdot L \in \mathcal{A}$. 

\smallskip

\noindent
Equality, for some $\lambda \in (0,1)$ and $p>1$, implies that $K$ and $L$ satisfy \eqref{e:BM functional F neg alpha} with equality for certain $\bar{\lambda}\in(0,1)$.
Moreover, if $\F$ is strictly decreasing and $K$ and $L$ are compact then $\conv K$ and $\conv L$ contain the origin and are dilatates of each other. If further $\F$ is strictly super-homogeneous of degree $1/\alpha$ then $\conv K = \conv L$.
\end{theorem}

\smallskip 

In the following, we present an example of an application of Theorem \ref{teorema F alpha negative}. To this aim, let $\mathcal{A}=\K^n_0$ be the family of all convex bodies containing the origin in their interior. Now, set $i\in\{0,1,\dots,n-1\}$, and define the functional $\mathcal{F}_i$ given by $\F_i(K) = \W_i(K^*)$ for any $K \in \K^n_0$. Let us verify that, indeed, $\mathcal{F}_i$ satisfies the hypotheses of Theorem \ref{teorema F alpha negative}:

\begin{enumerate}\itemsep5pt
\item Decreasing monotonicity: if \(K \subset L\), for \(K, L \in \K^n_0\), then \(K^* \supset L^*\) and thus
\begin{equation*}
\mathcal{F}_i(K) = \W_i(K^*) \geqslant \W_i(L^*) = \mathcal{F}_i(L),
\end{equation*}
which confirms that $\mathcal{F}_i$ is decreasing under set inclusion.
    
\item Super-homogeneity: for $r \geqslant 1$, scaling $K$ by $r$ yields $(rK)^* = {r}^{-1}K^*$. Hence,
\begin{equation*}
\mathcal{F}_i(rK) =\W_i\left(\frac{1}{r}K^*\right) = \left(\frac{1}{r}\right)^{n-i} \W_i(K^*) = r^{-(n-i)} \mathcal{F}_i(K).
\end{equation*}
Setting $\alpha = -1/(n-i)$ $\bigl(\text{and so } 1/\alpha = -(n-i)\bigr)$, this becomes
\begin{equation*}
\mathcal{F}_i(rK) = r^{1/\alpha} \mathcal{F}_i(K),
\end{equation*}
satisfying both \eqref{e: sub-homog} and \eqref{e: super-homoge} for all $r \geqslant 1$. 
    
\item Brunn-Minkowski type inequality: in \cite{F} Firey showed that, for $K, L \in \K^n_0$ and $\lambda \in (0,1)$,
\begin{equation}\label{e: BM Firey polar}
\W_i\Bigl(\bigl[(1-\lambda)K + \lambda L\bigr]^*\Bigr) \leqslant \left((1-\lambda)\W_i(K^*)^{-1/(n-i)} + \lambda\W_i(L^*)^{-1/(n-i)}\right)^{-(n-i)},
\end{equation}
with equality, for some  $\lambda\in(0, 1)$, if and only if $K$ and $L$ are dilatates. In other words, the latter inequality reads
\begin{equation*}
\mathcal{F}_i\bigl((1-\lambda)K + \lambda L\bigr) \leqslant \left((1-\lambda)\mathcal{F}_i(K)^{-1/(n-i)} + \lambda\mathcal{F}_i(L)^{-1/(n-i)}\right)^{-(n-i)},
\end{equation*}
which gives \eqref{e:BM functional F neg alpha} with $\alpha = -1/(n-i)$ and $C = 1$.
\end{enumerate}

Therefore, by Theorem \ref{teorema F alpha negative}, for any $p > 1$ and all $\lambda \in (0,1)$, the following $L_p$ Brunn-Minkowski inequality holds:
\begin{equation}\label{e:BM quermassintegrals Lp}
\W_i\Bigl(\bigl[(1-\lambda) \cdot K +_p \lambda \cdot L\bigr]^*\Bigr)^{-p/(n-i)} \geqslant \left((1-\lambda)\W_i(K^*)^{-p/(n-i)} + \lambda\W_i(L^*)^{-p/(n-i)}\right).
\end{equation}
Equality, for some $\lambda \in (0,1)$ and $p > 1$, implies that $K$ and $L$ satisfy \eqref{e: BM Firey polar} with equality, and thus they are dilatates. The sufficient condition of the equality case is straightforward.

This $L_p$ Brunn-Minkowski type inequality for polar bodies \eqref{e:BM quermassintegrals Lp}, jointly with its equality case, was previously obtained (following a different approach)  by Hernández Cifre and the second-named author in \cite{HY Polar}.

\section{L$_p$ Brunn-Minkowski inequalities for general absolutely continuous measures with radially decreasing densities}\label{s: L_p BM Gaussian}
\subsection{Background and preliminary results}
Since the original problem formulated by Gardner and Zvavitch in \cite{GZ} (cf. \eqref{e: Gaussian BM}), a vast array of results has been obtained in this line. Firstly, Gardner and Zvavitch proved such a Gaussian Brunn-Minkowski inequality for special families of sets, being some of these results later extended by Marsiglietti \cite{Ma} to the case of more general measures. In \cite{CLM}, Colesanti, Livshyts and Marsiglietti proved this inequality when both convex bodies $K$ and $L$ are small perturbations of the Euclidean ball. Subsequently, Livshyts, Marsiglietti, Nayar and Zvavitch in \cite{LMNZ} showed that the Brunn-Minkowski inequality is true for unconditional product measures with decreasing density and a pair of unconditional sets. This result was later generalized to the case of weakly unconditional sets by Ritoré and the second-named author \cite{RY}.
To introduce the precise statement of the latter result, first we recall the definition of weakly unconditional sets.
\begin{definition} \label{def: weakly sets}
We say that a set $A \subset \mathbb{R}^n$ is weakly unconditional if for any $(x_1, \ldots, x_n)\in A$ and all $(\epsilon_1, \ldots, \epsilon_n) \in \{0,1\}^n$ one has 
\begin{equation*}
 (\epsilon_1 x_1, \ldots, \epsilon_n x_n) \in A.
\end{equation*}
\end{definition}
We emphasize that this notion extends the well-known concept of unconditional sets: a subset $A \subset \mathbb{R}^n$ is said to be unconditional if for any $(x_1,\ldots, x_n)\in A$ and all $(\epsilon_1,\ldots,\epsilon_n)$ with $\epsilon_i \in [-1,1]$ one has $(\epsilon_1 x_1, \ldots, \epsilon_n x_n) \in A$.

We also need to recall the notion of a (strictly) radially decreasing function:
\begin{definition}
We say that a function $f: \R^n \longrightarrow \R_{\geqslant 0}$ is radially decreasing if
\begin{equation}\label{e: radially decreasing}
f(tx) \leqslant f(x)
\end{equation}
for all $x \in \R^n$ and any $t\geqslant 1$.
Furthermore, it is said to be strictly radially decreasing if the above inequality holds strictly whenever $t>1$.
\end{definition}
In the following, for a nonnegative function $f: \R^n \longrightarrow \R_{\geqslant 0}$, we will write $\mathrm{supp} f$ to denote the support of $f$, that is, the set $\{x\in\R^n: f(x)>0\}$. Moreover, we will say that $f$ is (strictly) radially decreasing on $\supp f$ if 
\eqref{e: radially decreasing} holds (strictly) whenever $x, tx\in\supp f$ (provided that $t>1$).

\smallskip

With these definitions, the above-mentioned result for weakly unconditional sets reads as follows:
\begin{thm}\label{t:Ritore Jesús measure BM}
Let $\nu=\nu_1 \times \cdots \times \nu_n$ be a product measure on $\R^n$ such that $\nu_i$ is the measure given by $\dlat \nu_i(x)=f_i(x) \dlat x$, where $f_i:\R \longrightarrow \R_{\geqslant 0}$ is a radially decreasing function, $i=1, \ldots, n$. 

Let $K, L \subset \mathbb{R}^n$ be weakly unconditional measurable sets with $\nu(K)\nu(L)>0$ such that $(1-\lambda)K+\lambda L$ is also measurable. Then, for all $\lambda \in (0,1)$,
\begin{equation*} 
\nu\bigl((1-\lambda)K+ \lambda L\bigr)^{1/n} \geqslant (1-\lambda ) \nu(K)^{1/n} + \lambda \nu(L)^{1/n}.
\end{equation*} 
\end{thm}
Moreover, Böröczky and Kalantzopoulos \cite{BK} showed that the dimensional Brunn-Minkowski ine\-quality \eqref{e: Gaussian BM} holds for convex bodies with symmetries with respect to $n$ independent hyperplanes:

\begin{thm}\label{t: Gaussian BM BK}
Let $H_1, \dots, H_n$ be (linear) hyperplanes with $H_1\cap \dots \cap H_n = \{ 0\}$. Let $K, L \in \K^n$ be convex bodies that are invariant under the orthogonal reflections through $H_1, \dots, H_n$. Then, for all $\lambda \in (0,1)$,
\begin{equation*}
\gamma_n\bigl((1-\lambda) K+ \lambda L\bigr)^{{1}/n} \geqslant (1-\lambda ) \gamma_n(K)^{{1}/n} + \lambda \gamma_n(L)^{1/n}.
\end{equation*} 
\end{thm}

In \cite{KL}, Kolesnikov and Livshyts followed a different approach and proved that \eqref{e: Gaussian BM} holds, with exponent $1/(2n)$ instead of $1/n$, when $K$ and $L$ are $0$-symmetric closed convex sets. Finally, by exploiting this latter approach, Eskenazis and Moschidis \cite{EM} obtained their celebrated result, Theorem \ref{t: EskMosch Gaussian BM}. 

\smallskip 

In this section, we explore $L_p$ Brunn-Minkowski inequalities for general absolutely continuous measures with radially decreasing density functions, being the Gaussian measure a particular case of them. 
To this end, we recall the analytical counterpart (for functions) of the Brunn-Minkowski inequality \eqref{e: BM}, the so-called Borell-Brascamp-Lieb inequality, originally proved in \cite{Bo} and \cite{BL} (which has as a particular case the well-known Prékopa-Leindler inequality). We refer the reader to \cite{G} for a detailed presentation of it, and we collect it here for the sake of completeness.
\begin{thm}[The Borell-Brascamp-Lieb inequality]\label{t: BBL}  Let $\lambda \in(0,1)$. Let $-1 / n \leqslant \beta \leqslant \infty$ and let $f, g, h: \mathbb{R}^n \longrightarrow \mathbb{R}_{\geqslant 0}$ be integrable functions with $\|f\|_1, \|g\|_1>0$ such that
\begin{equation}\label{e: condition BBL}
h\bigl((1-\lambda) x+\lambda y \bigr) \geqslant \Bigl((1-\lambda) f(x)^\beta+\lambda g(y)^\beta\Bigr)^{1 / \beta}
\end{equation}
for all $x, y \in \mathbb{R}^n$ with $f(x) g(y)>0$. Then
\begin{equation*}
\int_{\mathbb{R}^n} h(x) \mathrm{d} x \geqslant \Biggl[(1-\lambda)\left(\int_{\mathbb{R}^n} f(x) \;\dlat x\right)^\alpha+ \lambda \left(\int_{\mathbb{R}^n} g(x) \;\dlat x\right)^\alpha\Biggr]^{1/\alpha},
\end{equation*}
where $\alpha=\beta/(n\beta+1)$.
\end{thm}
Moreover, the special case of $\beta=\infty$ in
condition \eqref{e: condition BBL} of the previous result must be understood as its limit case, that is 
\begin{equation*}
h\bigl((1-\lambda) x+\lambda y \bigr) \geqslant \max \{f(x), g(y)\}.
\end{equation*}

Naturally connected to the Borell-Brascamp-Lieb inequality one finds the notion of $\beta$-concave function:
\begin{definition}
Let $\beta\in\R$, $\beta\neq0$. We say that a function $f:\R^n \longrightarrow \R_{\geqslant 0}$ is $\beta$-concave if 
\begin{equation*}
 f \bigl((1-\lambda)x + \lambda y \bigr) \geqslant \Bigl((1-\lambda) f(x)^\beta + \lambda f(y)^\beta\Bigr)^{1/\beta}
\end{equation*}
for all $x,y \in \R^n$ such that $f(x)f(y)>0$ and any $\lambda \in (0,1)$. Moreover, we say that $f$ is log-concave (sometimes also referred to as $0$-concave) if
\begin{equation*}
 f \bigl((1-\lambda)x + \lambda y \bigr) \geqslant f(x)^{1-\lambda}f(y)^{\lambda}
\end{equation*}
for all $x,y \in \R^n$ and any $\lambda \in (0,1)$.
\end{definition}
As a straightforward consequence of Theorem \ref{t: BBL}, we get the following Brunn-Minkowski inequality for measures with $\beta$-concave densities, for $\beta>0$ (see \cite[Corollary~11.2]{G}):

\begin{cor}\label{c: BM measures BBL}
Let $0< \beta$ and let $\nu$ be a measure on $\R^n$ given by $\dlat \nu (x) = f(x) \dlat x$, where $f$ is a $\beta$-concave function. Let $K, L\subset \R^n$ be measurable sets with $\nu(K)\nu(L)>0$ such that $(1-\lambda)K+\lambda L$ is also measurable. Then, for all $\lambda \in (0,1)$,
\begin{equation*} 
{\nu \bigl( (1-\lambda) K+ \lambda  L\bigr)^{\alpha} \geqslant (1-\lambda ) \nu(K)^{\alpha} + \lambda \nu(L)^{\alpha}},
\end{equation*}
where $\alpha=\beta/(n\beta+1)$.
\end{cor}

\begin{remark}\label{r: alpha smaller 1/n}
We note that, in the previous result, $\alpha$ belongs to $(0,1/n)$, being the limit case of $\alpha=1/n$ only possible (for arbitrary measurable sets $K, L\subset\R^n$ with $(1-\lambda)K+\lambda L$ measurable too) when $\nu$ is the restriction of the Lebesgue measure (up to a constant) to a convex set. Moreover, there is no absolutely continuous measure satisfying a Brunn-Minkowski inequality with degree of concavity $\alpha>1/n$ for any pair of measurable sets $K, L\subset\R^n$. We refer the reader to \cite{Bo, YN} and the references therein for more information on these questions.
\end{remark}

\smallskip

Another example of general absolutely continuous measures and sets for which some Brunn-Minkowski inequality (for the Minkowski addition) holds is the following elegant result by Livshyts, Marsiglietti, Nayar and Zvavitch in \cite{LMNZ}. We observe that, in the particular case of the Gaussian measure, this result holds true in arbitrary dimension because of the above-mentioned celebrated result by Eskenazis and Moschidis, collected in Theorem \ref{t: EskMosch Gaussian BM}.

\begin{thm}
Let $\nu$ be a measure on $\R^2$ given by $\dlat \nu (x) = f(x) \dlat x$, where $f$ is an even log-concave function. 
Let $K, L \subset \mathbb{R}^2$ be $0$-symmetric closed convex sets. Then, for all $\lambda \in (0,1)$,  
\begin{equation}\label{e: BM dim2 0-symmetric logc measures} 
\nu\bigl((1-\lambda) K+ \lambda  L\bigr)^{1/2} \geqslant (1-\lambda ) \nu(K)^{1/2} + \lambda \nu(L)^{1/2}.
\end{equation}
\end{thm}

We conclude this subsection by giving the following definition, which will be useful due to the assumption of Theorem \ref{teoremafunctionalFnonconvex} where a certain pair of sets satisfies a Brunn-Minkowski type inequality for the Minkowski addition:
\begin{definition}
Let $0<\alpha$ and let $\nu$ be an absolutely continuous measure on $\R^n$. Then, a pair $(K, L)$ of nonempty measurable sets with $\nu(K)\nu(L)>0$, such that any linear combination $\lambda_1 K+\lambda_2 L$, for $\lambda_1, \lambda_2\geqslant0$ with $\lambda_1+\lambda_2\leqslant 1$, is also measurable, is called $(\nu, \alpha)$-admissible if it satisfies that
\begin{equation*}\label{e:medidaB-M}
\nu\bigl((1-\lambda)K+ \lambda L\bigr)^{\alpha} \geqslant (1-\lambda ) \nu(K)^{\alpha} + \lambda \nu(L)^{\alpha}
\end{equation*}
for all $\lambda \in (0,1)$. 
\end{definition}

\subsection{Main results}
We start by showing the following consequence of Theorem \ref{teoremafunctionalFnonconvex}, from which we will derive different $L_p$ Brunn-Minkowski inequalities for general measures $\nu$ with radially decreasing densities, when dealing with a pair of $(\nu,\alpha)$-admissible sets, with $0<\alpha\leqslant 1/n$. We observe that here we are interested in Brunn-Minkowski inequalities with such possible values for the ``degree of concavity'', since unless one works with very particular families of subsets one cannot expect a stronger concavity for $\nu$ (see Remark \ref{r: alpha smaller 1/n}).

\begin{theorem}\label{t:LpBM_for_nu}
Let $0<\alpha \leqslant 1/n$ and let $\nu$ be a measure on $\R^n$ given by $\dlat \nu (x) = f(x) \dlat x$, where $f$ is a radially decreasing function with convex support $\supp f$. Let $(K, L)$ be a pair of $(\nu,\alpha)$-admissible sets and let $p > 1$. Then, for all $\lambda \in (0,1)$, 
\begin{equation} \label{e:LpBM_for_nu}
\nu\bigl((1-\lambda)\cdot K+_p \lambda \cdot L\bigr)^{{p}\alpha} \geqslant (1-\lambda ) \nu(K)^{{p}\alpha} + \lambda \nu(L)^{p\alpha}
\end{equation} 
whenever $(1-\lambda)\cdot K+_p\lambda\cdot L$ is measurable.

Equality, for some $\lambda \in (0,1)$ and $p>1$, when $K,L$ are also convex bodies contained in $\mathrm{supp} f$ and $f$ is further strictly radially decreasing on $\mathrm{supp} f$, holds if and only if $K=L$ and $0 \in K \cap L$.
\end{theorem}
\begin{proof}
Taking $\mathcal{A}$ as the family of all nonempty measurable sets in $\R^n$, we observe that $\mathcal{F}=\nu(\cdot)$
is increasing under set inclusion and sub-homogeneous of degree $1/\alpha$, since its density function $f$ is radially decreasing.  More precisely, doing the change of variables given by $x=ry$, jointly with the fact that $f$ is a radially decreasing function and $r \geqslant 1$, we get
\begin{equation}\label{e:homogeneity nu}
\nu(rA)=\int_{rA} f(x) \;\dlat x = r^n \int_{A} f(r y)\;\dlat y \leqslant r^n \int_{A} f(y)\;\dlat y =r^n \nu(A) \leqslant r^{1/\alpha} \nu(A),
\end{equation}
due to the fact that $\alpha \leqslant 1/n$.
Then, applying Theorem \ref{teoremafunctionalFnonconvex} to $\F$ and the pair $(K, L)$ (see also Remark \ref{r: A not necessarily closed under linear comb}) we obtain \eqref{e:LpBM_for_nu}. 

For the equality case, first we observe that since $\supp f$ is a convex set, which further contains the origin without loss of generality, if $K,L$ are convex sets that are contained in $\supp f$ then (following the notation $t=t(\mu)$, $s=s(\mu)$, for a fixed $\lambda\in(0,1)$, of the proof of Theorem \ref{teoremafunctionalFnonconvex}) we have that the sets 
\[
(1-\lambda)\cdot K+_p \lambda\cdot L = \bigcup_{\mu \in [0,1]}\bigl( t(\mu)K+s(\mu)L\bigr), \quad tK+sL \quad \text{and} \quad \frac{t}{t+s}K+\frac{s}{t+s}L
\]
are contained in $\supp f$, because $t+s\leqslant 1$.
Thus, the necessity of the equality case is derived from (the proof of) Theorem \ref{teoremafunctionalFnonconvex}, by taking into account that \eqref{e:homogeneity nu} holds strictly whenever $rA\subset\supp f$ and $r>1$, jointly with the compactness of $K$ and $L$ and the fact that $\nu$ is strictly increasing when dealing with convex sets with nonempty interior.
Finally, when $K,L\in\K^n$ are convex bodies containing the origin, we clearly have that if $K=L$ then $(1-\lambda)\cdot K +_p\lambda\cdot L=K=L$ (cf. \eqref{e: p-sum Firey}) and hence \eqref{e:LpBM_for_nu} holds with equality.
\end{proof}

Regarding the original problem by Gardner and Zvavitch in \cite{GZ} (cf. \eqref{e: Gaussian BM}), and aiming to obtain an $L_p$ version of such an inequality, we obtain the following consequence of our previous result. This is due to the fact that the Gaussian measure is strictly sub-homogeneous, because its density is strictly radially decreasing.

\begin{theorem}\label{t:gamma-admissible}
Let $(K, L)$ be a pair of $(\gamma_n,1/n)$-admissible sets and $p > 1$. Then, for all $\lambda \in (0,1)$, 
\begin{equation*}\label{e:LpBM_for_gamma}
\gamma_n\bigl((1-\lambda)\cdot K+_p \lambda \cdot L\bigr)^{{p}/n} \geqslant (1-\lambda ) \gamma_n(K)^{{p}/n} + \lambda \gamma_n(L)^{p/n}
\end{equation*} 
whenever $(1-\lambda)\cdot K+_p\lambda\cdot L$ is measurable.

Equality, for some $\lambda \in (0,1)$ and $p>1$, when $K,L$ are also convex bodies, holds if and only if $K=L$ and $0 \in K \cap L$.
\end{theorem}

In the following we collect some examples of application of Theorem \ref{t:LpBM_for_nu} by using some known results for the Minkowski addition. 
First, we observe that Theorem \ref{t:LpBM_for_nu} combined with Theorem \ref{t:Ritore Jesús measure BM} leads us to the next result, which has been recently established (without the equality conditions) by Roysdon and Xing in \cite{RX}, by using a different approach.
\begin{theorem}\label{t: Gaussian Lp BM weakly}
Let $\nu=\nu_1 \times \cdots \times \nu_n$ be a product measure on $\R^n$ such that $\nu_i$ is the measure given by $\dlat \nu_i(x)=f_i(x) \dlat x$, where $f_i:\R \longrightarrow \R_{\geqslant 0}$ is a radially decreasing function with convex support $\supp f_i$, $i=1, \ldots, n$. 

Let $K, L \subset \mathbb{R}^n$ be weakly unconditional measurable sets with $\nu(K)\nu(L)>0$ such that any linear combination $\lambda_1 K+\lambda_2 L$, for $\lambda_1, \lambda_2\geqslant0$ with $\lambda_1+\lambda_2\leqslant 1$, is also measurable and let $p > 1$. Then, for all $\lambda \in (0,1)$,
\begin{equation} \label{e: nu Lp BM weakly}
\nu \bigl( (1-\lambda)\cdot K+_p \lambda \cdot L\bigr)^{p/n} \geqslant (1-\lambda ) \nu(K)^{p/n} + \lambda \nu(L)^{p/n}
\end{equation} 
whenever $(1-\lambda)\cdot K+_p\lambda\cdot L$ is measurable.

Equality, for some $\lambda \in (0,1)$ and $p>1$, when $K,L$ are also convex bodies contained in $\mathrm{supp} f$ and $f$ is further strictly radially decreasing on $\mathrm{supp} f$, holds if and only if $K=L$.
\end{theorem}

We observe that the Gaussian measure is a particular case of the above-considered product measures with radially decreasing densities. Thus, \eqref{e: nu Lp BM weakly} holds for $\gamma_n(\cdot)$ when dealing with weakly unconditional measurable sets. Another example of family of sets for which such a Gaussian Brunn-Minkowski holds is that of convex bodies with symmetries with respect to some hyperplanes. This is the content of the following result, which is a consequence of Theorem \ref{t:LpBM_for_nu} combined with Theorem \ref{t: Gaussian BM BK}.
\begin{theorem}\label{t: Gaussian Lp BM BK}
Let $H_1, \dots, H_n$ be (linear) hyperplanes with $H_1\cap \dots \cap H_n = \{ 0\}$. Let $K, L \in \K^n$ be convex bodies that are invariant under the orthogonal reflections through $H_1, \dots, H_n$ and let $p> 1$. Then, for all $\lambda \in (0,1)$,
\begin{equation*}
\gamma_n\bigl((1-\lambda)\cdot K+_p \lambda \cdot L\bigr)^{{p}/n} \geqslant (1-\lambda ) \gamma_n(K)^{{p}/n} + \lambda \gamma_n(L)^{p/n}.
\end{equation*} 
Equality, for some $\lambda \in (0,1)$ and $p>1$, holds if and only if $K=L$ and $0 \in K \cap L$.
\end{theorem}

Now, taking into account the application of the Borell-Brascamp-Lieb inequality collected in Corollary \ref{c: BM measures BBL} and Theorem \ref{t:LpBM_for_nu}, we arrive at the following result for measures with $\beta$-concave density functions, for $\beta>0$, with supremum at the origin. We notice that the radial decreasing monotonicity of the density function $f$ is guaranteed by these assumptions:
\begin{equation}\label{e: rad decreasing when max at the origin}
f(x)=f\left[\frac{1}{r}(rx)+\left(1-\frac{1}{r}\right)0\right]
\geqslant
\left[\frac{1}{r}f(rx)^{\beta}+\left(1-\frac{1}{r}\right)f(0)^{\beta}\right]^{1/\beta}
\geqslant f(rx)
\end{equation}
for all $x, rx\in\supp f$, whenever $r\geqslant 1$.
We point out that this result was recently shown by Roysdon and Xing in \cite{RX} (without the equality conditions) with a different approach.

\begin{theorem}\label{t: Lp BM measures}
Let $0<\alpha<1/n$ and let $\nu$ be a measure on $\R^n$ given by $\dlat \nu (x) = f(x) \dlat x$, where $f$ is a $\beta$-concave function with $f(0)=\sup_{x\in\R^n}f(x)$, for $\beta=\alpha/(1-n\alpha)$. 

Let $K, L\subset \R^n$ be measurable sets with $\nu(K)\nu(L)>0$ such that any linear combination $\lambda_1 K+\lambda_2 L$, for $\lambda_1, \lambda_2\geqslant0$ with $\lambda_1+\lambda_2\leqslant 1$, is also measurable and let $p > 1$. Then, for all $\lambda \in (0,1)$,
\begin{equation*} 
{\nu \bigl( (1-\lambda)\cdot K+_p \lambda \cdot L\bigr)^{p\alpha} \geqslant (1-\lambda ) \nu(K)^{p\alpha} + \lambda \nu(L)^{p\alpha}}
\end{equation*}
whenever $(1-\lambda)\cdot K+_p\lambda\cdot L$ is measurable.

Equality, for some $\lambda \in (0,1)$ and $p>1$, when $K,L$ are also convex bodies contained in $\mathrm{supp} f$ and $f$ is further  strictly radially decreasing on $\mathrm{supp} f$, holds if and only if $K=L$ and $0 \in K \cap L$.
\end{theorem}

Again, the corresponding $L_p$ version of \eqref{e: BM dim2 0-symmetric logc measures}  may be obtained as a consequence of Theorem \ref{t:LpBM_for_nu}. The symmetry of the density function jointly with the log-concavity implies that its supremum is attained at the origin and then, using the same argument as the one in
\eqref{e: rad decreasing when max at the origin}, we conclude that it is radially decreasing. So, we get:

\begin{theorem}
Let $\nu$ be a measure on $\R^2$ given by $\dlat \nu (x) = f(x) \dlat x$, where $f$ is an even log-concave function. 
Let $K, L \subset \mathbb{R}^2$ be $0$-symmetric closed convex sets with $\nu(K)\nu(L)>0$ and $p > 1$. Then, for all $\lambda \in (0,1)$,
\begin{equation*} 
{\nu \bigl( (1-\lambda)\cdot K+_p \lambda \cdot L\bigr)^{p/2} \geqslant (1-\lambda ) \nu(K)^{p/2} + \lambda \nu(L)^{p/2}}.
\end{equation*}
Equality, for some $\lambda \in (0,1)$ and $p>1$, when $K,L$ are also compact sets contained in $\mathrm{supp} f$ and $f$ is further strictly radially decreasing on $\mathrm{supp} f$, holds if and only if $K=L$.
\end{theorem}

\medskip 

We conclude this subsection by noticing that by combining Theorem \ref{teoremafunctionalFnonconvex} with Theorem \ref{t: EskMosch Gaussian BM} we directly obtain Theorem \ref{t:Gaussian measure symmetric convex Lp BM}. 
We point out that the equality condition, for $p>1$, follows directly from the characterization of the equality case given in Theorem \ref{t: EskMosch Gaussian BM}. However, in the following subsection we provide a proof of the equality case for $p>1$ independently of the characterization for $p=1$.

\subsection{An alternative proof of (the equality case of) Theorem \ref{t:Gaussian measure symmetric convex Lp BM}}

We notice that \eqref{e:Gaussian measure symmetric convex Lp BM}  immediately follows from Theorem \ref{t:gamma-admissible} jointly with Theorem \ref{t: EskMosch Gaussian BM}.
In the following we present an argument that allows us to derive further the equality case without using that of Theorem \ref{t: EskMosch Gaussian BM} (i.e., without using the characterization of the equality condition for the inequality for the Minkowski addition).

So, suppose that \eqref{e:Gaussian measure symmetric convex Lp BM} holds with equality. First, we observe that if $K$ and $L$ are bounded then $K=L$, by the equality case of Theorem \ref{t:gamma-admissible}, and we are done. Hence, assume that at least one of the sets, say $K$, is unbounded. Since $K$ is $0$-symmetric, there exists an $(n-k)$-plane $H$, with $0\leqslant k\leqslant n-1$, such that 
\begin{equation}\label{e:K=K|H^{perp}+H}
K=K|H^{\perp}+H,
\end{equation}
where $K|H^{\perp}$ is a bounded set. Indeed, since $K$ is convex and unbounded, $K$ contains a ray $R=\{x+ry: r\geqslant0\}$, for some $x,y\in\R^n$ (see, e.g. \cite[Exercise 13, Section 1.1]{HW}), and then $\conv(R\cup\{0\})\subset K$. The latter implies that the one-dimensional vector subspace $H_1=\{sy: s\in\R\}$ is contained in $K$, and moreover $K=K|H^{\perp}_1 + H_1$, because of the central symmetry of $K$ (and taking into account that $K$ is closed). Then, if $K|H^{\perp}_1$ is bounded we are done; otherwise, we proceed analogously with $K|H^{\perp}_1$, obtaining that there exists a one-dimensional vector subspace $H_2$ such that 
\[
K|H^{\perp}_1=\bigl(K|H^{\perp}_1\bigr)\big|H^{\perp}_2 + H_2.
\]
By repeating this process finitely many times, we get \eqref{e:K=K|H^{perp}+H}.

Notice that if $k=0$ then $K=\R^n$, and so having equality in \eqref{e:Gaussian measure symmetric convex Lp BM} implies that $\gamma_n(L)$ must be one. Hence we deduce that $L=\R^n$ because $L$ is convex, and thus $K=L$. Therefore, we may assume that $1\leqslant k\leqslant n-1$.

On the one hand, applying Theorem \ref{t:Gaussian measure symmetric convex Lp BM} with $K|H^{\perp}$ and $L|H^{\perp}$, since such projections are $0$-symmetric closed convex sets with nonempty interior (in $\R^k$), we have that 
\begin{equation}\label{e: BM gammak proy}
\gamma_k\Bigl((1-\lambda)\cdot\bigl(K|H^{\perp}\bigr)+_p\lambda\cdot\bigl(L|H^{\perp}\bigr)\Bigr)
\geqslant \Bigl((1-\lambda)\gamma_k\bigl(K|H^{\perp}\bigr)^{p/k}+
\lambda\gamma_k\bigl(L|H^{\perp}\bigr)^{p/k}\Bigr)^{k/p}.
\end{equation}
Moreover, due to the fact that $\gamma_n$ is a rotationally invariant product measure, $\gamma_n(K)=\gamma_k\bigl(K|H^{\perp}\bigr)$. Thus, 
\begin{equation}\label{e: L subset proj}
\begin{split}
\Bigl((1-\lambda)\gamma_k\bigl(K|H^{\perp}\bigr)^{p/k}+
\lambda\gamma_k\bigl(L|H^{\perp}\bigr)^{p/k}\Bigr)^{k/p}
&=\Bigl((1-\lambda)\gamma_n(K)^{p/k}+
\lambda\gamma_n\bigl(L|H^{\perp}+H\bigr)^{p/k}\Bigr)^{k/p}\\
&\geqslant\Bigl((1-\lambda)\gamma_n(K)^{p/k}+
\lambda\gamma_n(L)^{p/k}\Bigr)^{k/p}, 
\end{split}
\end{equation}
where the inequality follows from $L\subset L|H^{\perp}+H$.
Now, by the monotonicity of means (see \cite[Theorem $16$]{HLP}) 
jointly with the equality case of \eqref{e:Gaussian measure symmetric convex Lp BM}, we get
\begin{equation}\label{e: monotonicity means gamma}
\begin{split}
\Bigl((1-\lambda)\gamma_n(K)^{p/k}+
\lambda\gamma_n(L)^{p/k}\Bigr)^{k/p}
&\geqslant
\Bigl((1-\lambda)\gamma_n(K)^{p/n}+
\lambda\gamma_n(L)^{p/n}\Bigr)^{n/p}\\
&=\gamma_n\bigl((1-\lambda)\cdot K+_p \lambda\cdot L\bigr).
\end{split}
\end{equation}

On the other hand, denoting by
\begin{equation*}
M=\left\{(1-\mu)^{1/q}(1-\lambda)^{1/p}x+\mu^{1/q}\lambda^{1/p}y: x\in K|H^{\perp}, y\in L|H^{\perp}, \, \mu\in[0,1)\right\},
\end{equation*}
we have that
\begin{equation}\label{e: M+H}
M+H\subset (1-\lambda)\cdot K+_p \lambda\cdot L.
\end{equation}
To see this, let $z\in M+H$. By the definition of $M$, there exist $\mu\in[0,1)$, $x_1\in K|H^{\perp}$, $y_1\in L|H^{\perp}$ and $z_2\in H$
such that 
\[
z=(1-\mu)^{1/q}(1-\lambda)^{1/p}x_1+\mu^{1/q}\lambda^{1/p}y_1 + z_2.
\]
Furthermore, since $y_1\in L|H^{\perp}$, there exists $y_2\in H$ such that $y_1+y_2\in L$. Then, setting $\alpha=(1-\mu)^{1/q}(1-\lambda)^{1/p}>0$, we clearly have that 
\begin{equation*}
z=(1-\mu)^{1/q}(1-\lambda)^{1/p}\bigl(x_1+(1/\alpha)(z_2-\mu^{1/q}\lambda^{1/p}y_2)\bigr)+\mu^{1/q}\lambda^{1/p}(y_1+y_2)\in (1-\lambda)\cdot K+_p \lambda\cdot L,
\end{equation*}
due to the fact that $K=K|H^{\perp}+H$.

Observe also that, since $\lambda\cdot(L|H^{\perp})$ is a convex set containing the origin, we get
\begin{equation}\label{e: inclusions int proj}
\relinter \Bigl( \lambda \cdot \bigl(L|H^{\perp}\bigr) \Bigr)  \subset \bigcup_{\mu \in [0,1)} \Bigl( \mu^{1/q} \lambda \cdot \bigl(L|H^{\perp} \bigr) \Bigr)\subset M.
\end{equation}
Hence, since $M$ is clearly contained in $(1-\lambda)\cdot\bigl(K|H^{\perp}\bigr)+_p\lambda\cdot\bigl(L|H^{\perp}\bigr)$, and taking in account \eqref{e: inclusions int proj}, it follows that $M$ differs from $(1-\lambda)\cdot\bigl(K|H^{\perp}\bigr)+_p\lambda\cdot\bigl(L|H^{\perp}\bigr)$ in a set of $\gamma_k$-measure zero, because $\lambda \cdot \bigl(L|H^{\perp}\bigr)$ is convex.
Therefore, we conclude that 
\begin{equation*}
\gamma_k(M)
=\gamma_k\Bigl((1-\lambda)\cdot\bigl(K|H^{\perp}\bigr)+_p\lambda\cdot\bigl(L|H^{\perp}\bigr)\Bigr).
\end{equation*}
So, \eqref{e: M+H} gives
\begin{equation}\label{e: gamma k p proj sum}
\gamma_n\bigl((1-\lambda)\cdot K+_p \lambda\cdot L\bigr)\geqslant\gamma_n(M+H)=
\gamma_k(M)
=\gamma_k\Bigl((1-\lambda)\cdot\bigl(K|H^{\perp}\bigr)+_p\lambda\cdot\bigl(L|H^{\perp}\bigr)\Bigr).
\end{equation}

Now, joining \eqref{e: BM gammak proy}, \eqref{e: L subset proj}, \eqref{e: monotonicity means gamma} and \eqref{e: gamma k p proj sum}, there must be equality in all these inequalities.
Then, assuming that $L|H^{\perp}$ is bounded, equality in \eqref{e: BM gammak proy} 
yields $K|H^{\perp}=L|H^{\perp}$, due to the case, shown at the beginning of the proof, where both sets are bounded. Furthermore, equality in \eqref{e: L subset proj} gives $L=L|H^{\perp}+H$ and, therefore, we have that $K=L$. Finally, the case where $L|H^{\perp}$ is unbounded, is obtained similarly, proceeding with $L|H^{\perp}$ as above for $K$. This finishes the proof.

\begin{remark}
Let $(K,L)$ be a pair of $(\gamma_n, 1/n)$-admissible sets satisfying with equality the $L_p$ Gau\-ssian Brunn-Minkowski inequality \eqref{e:Gaussian measure symmetric convex Lp BM}. If there exists an $(n-k)$-plane $H$ such that $K=K|H^{\perp}+H$ and $L=L|H^{\perp}+H$, being $K|H^{\perp}$ and $K|H^{\perp}$ bounded sets, then $K=L$, provided that $(K|H^{\perp},L|H^{\perp})$ is a pair of $(\gamma_k, 1/k)$-admissible sets.
\end{remark}


\section{L$_p$ Brunn-Minkowski inequalities for the Wills functional}\label{s: BM Wills functional}

In this section, we derive some consequences of our main result, Theorem \ref{teoremafunctionalFnonconvex}, for the (generalized) Wills functional. Along this section, unless we say the opposite, all the sets involved ($K, L, E\subset\R^n$) will be convex bodies.

Firstly, we obtain the corresponding $L_p$ Brunn-Minkowski inequality for the Wills functional 
(see Theorem \ref{t:classic Wills functional Lp BM}) as a consequence of Theorem \ref{teoremafunctionalFnonconvex}, due to the fact that the Wills functional satisfies the assumptions therein. Indeed, on the one hand, $\w(\cdot)$ is increasing under set inclusion, because of the monotonicity of the intrinsic volumes, i.e., if $K \subset L$, then $\w(K) \leqslant \w(L)$.
Moreover, it is sub-homogeneous of degree $n$, namely, 
\begin{equation*}\label{e:Willshomogeneous}
    \w(r K) \leqslant r^n \w(K)
\end{equation*}
for all $r \geqslant 1$. Although the proof follows an approach similar to the one presented earlier for an absolutely continuous measure with a radially decreasing density function (discussed in Section \ref{s: L_p BM Gaussian}), we include the full argument for the sake of completeness. Hence, to see the sub-homogeneity of the Wills functional, we will make use of the following integral representations showed by Hadwiger in \cite{H}:  
\begin{equation}\label{e: integralWillsHadwiger}
 \w(K)=\int_{\R^n} \e^{-\pi d(x,K)^2}\;\mathrm{d}x= 2\pi \int_{0}^{\infty} \vol(K+tB_n)t\e^{-\pi t^2} \;\mathrm{d}t,
\end{equation} 
where $d(x, K)=\min_{y\in K} \|x-y\|$ is the Euclidean distance from $x$ to $K$.
Then, using that the distance is homogeneous of degree one, we get
\begin{equation}\label{e: proving subh Wills I}
    {\w(r K)=\int_{\mathbb{R}^n} \e^{-\pi {d\left({x}, rK\right)^2 }}\dlat x}=\int_{\mathbb{R}^n} \e^{-\pi r^2  {d \left(\frac{x}{r}, K\right)^2} }\dlat x
\end{equation}
for any $r\geqslant 1$.
Now, doing the change of variables given by $y = x / r$ and using that the negative exponential is a decreasing function jointly with the fact that $r \geqslant 1$, we obtain
\begin{equation}\label{e: proving subh Wills II}
\int_{\mathbb{R}^n} \e^{-\pi r^2  {d \left(\frac{x}{r}, K\right)^2} }\dlat x =
r^n \int_{\mathbb{R}^n} \e^{-\pi r^2  {d \left(y, K\right)^2} }\dlat y
\leqslant r^n \int_{\mathbb{R}^n} \e^{-\pi  {d \left(y, K\right)^2} }\dlat y=r^n\w(K),
\end{equation}
and thus $\w(\cdot)$ is indeed sub-homogeneous of degree $n$.

On the other hand,  Theorem \ref{t: BM Wills} provides us with the corresponding Brunn-Minkowski type ine\-quality (with exponent $\alpha=1/n$) for any pair of convex bodies $K$ and $L$.

Therefore, altogether, applying Theorem \ref{teoremafunctionalFnonconvex} we derive Theorem \ref{t:classic Wills functional Lp BM}.

\smallskip

Now, we show that this result can be stated in a more general setting, that is, for the so-called generalized Wills functional, introduced by Kampf in \cite{K}, and previously studied for different questions (see e.g. \cite{AHY, HY, K} and the references therein). In this regard, Kampf \cite{K} established certain generalizations of the above integral formulas \eqref{e: integralWillsHadwiger}, by considering the `distance' $d_E(x,K)$, between $x \in \mathbb{R}^n $ and $K$, relative to a convex body $E$ with {$0 \in \inter E$}. More precisely, for
\begin{equation*}
    d_E(x, K)=\min_{y\in K} \|x-y\|_E=\min\{t \geqslant 0 \;:\;x \in K+tE\},
\end{equation*}  
he showed that
\begin{equation*}\label{e:integralWillsK}
 \int_{\R^n} \e^{-\pi d_E(x,K)^2}\;\mathrm{d}x= 2\pi \int_{0}^{\infty} \vol(K+tE)t\e^{-\pi t^2} \;\mathrm{d}t.
\end{equation*} 
It is easy to see that the latter equals $\sum_{i=0}^n \V_{i}(K)$ when $E$ is the unit ball $B_n$. 

In this context, for $\lambda \geqslant 0$, the well-known (relative) Steiner formula (cf. \eqref{eq:Steiner formula}) says that 
\begin{equation}\label{eq:Steiner relative formula}
    \vol(K+\lambda E)= \sum_{i=0}^n\binom{n}{i} \W_{i}(K;E) \lambda^i.
\end{equation}
The coefficients $\W_i(K;E)$, for convex bodies $K, E \in \K^n$, are the relative quermassintegrals of $K$ with respect to $E$.

It can be observed that replacing the density $t\e^{-\pi t^2}$ with a different function $G(t)$, appropria\-tely associated to a measure $\nu$ defined on the nonnegative real line $[0, \infty)$, yields a more general functional, which will be denoted by $\w^{\,\nu}(K;E)$. To be more precise, we have that
\begin{equation}\label{e: def generalized Wills funct}
\w^{\,\nu}(K;E)=\int_{\R^n} G \bigl(d_E(x, K) \bigr)\;\mathrm{d}x,  
\end{equation}
where $G(t)=\nu\bigl([t, \infty) \bigr)$ for any $t\in[0,\infty)$.
Now, exploiting the idea of the proof given by Kampf \cite{K} (see also \cite{AHY} and \cite{HY}), we have the following: if $\nu$ is given by $\mathrm{d}\nu(t)=\phi(t)\mathrm{d}t$, where $\phi: \mathbb{R}_{\geqslant 0} \longrightarrow  \mathbb{R}_{\geqslant 0}$ is a nonnegative measurable function with $\phi(t)>0$ for all $t \in [0, \infty)$, then $G(t)=\nu\bigl([t, \infty) \bigr)$ is strictly decreasing and moreover $G(t) >0 $ for all $t \in [0, \infty)$.
Hence, we notice that we may write 
\begin{equation*}
G(t)=\e^{-u(t)}
\end{equation*}
for some function $u:\R_{\geqslant 0} \longrightarrow \R$, namely, we can define $u$ as the function given by $u(t)=-\log G(t)$. So, 
\eqref{e: def generalized Wills funct} reads 
\begin{equation}\label{e: def generalized Wills funct exponential}
\w^{\,\nu}(K;E)=\int_{\R^n} \e^{-u\bigl(d_E(x, K)\bigr)}\;\mathrm{d}x 
\end{equation}
whenever the density of $\nu$ is strictly positive on $[0,\infty)$, and then from now on $\w^{\,\nu}(\cdot\,;E)$ will be denoted by $\w_{u}(\cdot\,;E)$ for the sake of consistency of the notation. 

Notice that the most general functional given by the integral expression \eqref{e: def generalized Wills funct} (for a certain function $G$) seems actually to be of the form of \eqref{e: def generalized Wills funct exponential}, since the decreasing monotonicity of $G$ is needed to assure one of the assumptions of our main result (Theorem \ref{teoremafunctionalFnonconvex}).

Furthermore, we point out that we clearly have
\begin{equation*}
\bigl\{x\in\R^n\,:\, u\bigl(d_E(x,K)\bigr)\leqslant s\bigr\}=\bigl\{x\in \R^n\,:\, d_E(x,K)\leqslant u^{-1}(s)\bigr\}=K + u^{-1}(s)E.
\end{equation*}
Thus, we get that 
\begin{equation}\label{e: generalized Wills weights}
\begin{split}
\w_{u}(K;E)&=\int_{\R^n} \int_{u\bigl(d_E(x,K)\bigr)}^{\infty} \e^{-s}\dlat s\; \dlat x 
= \int_{u(0)}^{\infty} \int_{\R^n}  \chi_{\{y\in\R^n\,:\, u(d_E(y,K))\leqslant s\}}(x) \, \dlat x \, \e^{-s}\dlat s\\
&= \int_{u(0)}^{\infty} \vol\bigl(K+u^{-1}(s)E\bigr) \; \e^{-s}\dlat s
= \sum_{i=0}^n \binom{n}{i}\W_{i}(K;E) \int_{u(0)}^{\infty} u^{-1}(s)^i \; \e^{-s}\dlat s,
\end{split}
\end{equation}
where we have used Fubini's theorem and the Steiner formula \eqref{eq:Steiner relative formula}, jointly with the fact that we may assure that there exists the inverse of $u$ (on its image) because  $u$ is strictly increasing.

In this setting, in \cite[Theorem 4.1]{AHY} the following Brunn-Minkowski type inequality for the genera\-lized Wills functional $\w_u(\,\cdot\, ;E)$ was shown. Before recalling its statement, we point out that there is a small typo in \cite[(4.2)]{AHY}. This is because, in its proof, one considers the measure on $[u(0),\infty)$ with density function $s\mapsto\e^{-s}$ (which is not, in principle, a probability measure), and then the nece\-ssary normalization $\e^{-s}\e^{u(0)}$ of the density, when applying Jensen's inequality, is missing. And this yields the (erroneous) multiplicative factor $\e^{-(n-1)u(0)/n}$ in the statement of the correspon\-ding Brunn-Minkowski inequality (see \cite[(4.2)]{AHY}). Hence, the (right) statement (when assuming only that $u$ is strictly increasing, which is actually the sole assumption needed there) of this result reads as follows:

\begin{thm}\label{t: generalized Wills BM}
Let $K, L, E \in \K^n$ be convex bodies with $0 \in \mathrm{int}\, E$ and let $u:\R_{\geqslant 0} \longrightarrow \R$ be a strictly increasing function. Then, for all $\lambda \in (0,1)$,
\begin{equation*}
    \w_u \bigl((1-\lambda)K + \lambda L; E \bigr)^{1/n} \geqslant \frac{1}{(n!)^{1/n}} \bigl( (1-\lambda) \w_u(K;E)^{1/n}+ \lambda \w_u(L;E)^{1/n}\bigr).
\end{equation*}
\end{thm}
We conclude the paper by noting that the generalized Wills functional $\w_u(\,\cdot\, ;E)$ (for a fixed convex body $E$ with nonempty interior), given by \eqref{e: def generalized Wills funct}, satisfies the conditions of Theorem \ref{teoremafunctionalFnonconvex}. Firstly, it is clear that such a functional is increasing. Indeed, if $K\subset L$ then $d_E(x,K)\geqslant d_E(x,L)$ for any $x\in \R^n$, and thus $G\bigl(d_E(x,K)\bigr) \leqslant G\bigl(d_E(x,L)\bigr)$, because of the monotonicity of $G$, from which we get $\w_u(K;E)\leqslant \w_u(L;E)$. Moreover, the generalized  Wills functional is also sub-homogeneous, i.e.,
\begin{equation*}\label{e:Willsghomogeneous}
\w_u(r K;E) \leqslant r^n \w_u(K;E)
\end{equation*}
if $r \geqslant 1$. The proof of this fact may be done with the same steps to those followed in \eqref{e: proving subh Wills I} and \eqref{e: proving subh Wills II}, by taking into account that $G$ is a (strictly) decreasing function. Then, combining Theorem \ref{teoremafunctionalFnonconvex} with Theorem \ref{t: generalized Wills BM}, we conclude the following $L_p$ Brunn-Minkowski inequality for the generalized Wills functional:

\begin{theorem}
Let $K,L, E \in \K^n $ be convex bodies with $0 \in \mathrm{int}\,E$, let $u:\R_{\geqslant 0} \longrightarrow \R$ be a strictly increasing function and $p> 1$. Then, for all $\lambda \in (0,1)$,
\begin{equation}\label{e:LpBMforwills generalized}
{\w_u \bigl((1-\lambda)\cdot K+_p \lambda \cdot L; E\bigr)^{{p}/{n}} \geqslant \frac{1}{(n!)^{{p}/{n}}} \left((1-\lambda ) \w_u(K;E)^{{p}/{n}} + \lambda \w_u(L;E)^{{p}/{n}}\right)}.
\end{equation}
\end{theorem}

\begin{remark}
Observe that the above inequality for the (generalized) Wills functional cannot hold with equality. Indeed, $\w_u(\,\cdot\, ;E)$ is strictly increasing, due to \eqref{e: generalized Wills weights} jointly with the monotonicity properties of the (relative) quermassintegrals under set inclusion (cf. \eqref{eq:Steiner relative formula}), and it is further strictly sub-homogeneous, because $G$ is strictly decreasing (cf. \eqref{e: proving subh Wills I} and \eqref{e: proving subh Wills II}). Therefore, from the equality case of Theorem \ref{teoremafunctionalFnonconvex}, equality in \eqref{e:LpBMforwills generalized} would imply that $K=L$. However, when $K=L$, one clearly does not have equality there by the presence of the constant ${1}/{(n!)^{p/n}}$. 
\end{remark}

\vspace{6mm}

\noindent {\it Acknowledgements.}
We would like to thank Dylan Langharst for bringing to our attention some interesting observations after the first version of our manuscript.

\end{document}